\documentclass[reqno,centertags, 12pt]{amsart}
\usepackage{amsmath,amsthm,amscd,amssymb}
\usepackage{latexsym,verbatim}

\newcommand{\bbR}{{\mathbb{R}}}

\newcommand{\bbP}{{\mathbb{P}}}
\newcommand{\bbE}{{\mathbb{E}}}
\newcommand{\bbZ}{{\mathbb{Z}}}

\newcommand{\bbT}{{\mathbb{T}}}

\newcommand{\lb}{\label}

\newcommand{\beq}{\begin{equation}}
\newcommand{\eeq}{\end{equation}}
\newcommand{\ba}{\begin{align}}
\newcommand{\ea}{\end{align}}
\newcommand{\eps}{\varepsilon}
\newcommand{\del}{\delta}
\newcommand{\tht}{\theta}
\newcommand{\ka}{\kappa}
\newcommand{\al}{\alpha}

\newcommand{\ga}{\gamma}

\newcommand{\til}{\tilde}
\newcommand{\tilth}{\til\tht}
\newcommand{\tilT}{\til T}

\newcounter{smalllist}
\newenvironment{SL}{\begin{list}{{\rm\roman{smalllist})}}{%
\setlength{\topsep}{0mm}\setlength{\parsep}{0mm}\setlength{\itemsep}{0mm}%
\setlength{\labelwidth}{2em}\setlength{\leftmargin}{2em}\usecounter{smalllist}%
}}{\end{list}}

\DeclareMathOperator*{\divg}{div}

\allowdisplaybreaks
\numberwithin{equation}{section}

\newtheorem{theorem}{Theorem}[section]

\newtheorem{lemma}[theorem]{Lemma}
\newtheorem{corollary}[theorem]{Corollary}
\theoremstyle{definition}
\newtheorem{definition}[theorem]{Definition}

\theoremstyle{remark}

\begin{document}
\title[Quenching and Propagation of Combustion]{Quenching and Propagation of Combustion Without Ignition Temperature
Cutoff}

\author{Andrej Zlato\v s}

\address{ Department of Mathematics \\ University of
Wisconsin \\ Madison, WI 53706, USA \\ Email: \tt
zlatos@math.wisc.edu}


\maketitle

\begin{abstract}
We study a reaction-diffusion equation in the cylinder $\Omega =
\mathbb{R}\times\mathbb{T}^m$, with combustion-type reaction term
without ignition temperature cutoff, and in the presence of a
periodic flow. We show that if the reaction function decays as a
power of $T$ larger than three as $T\to 0$ and the initial datum
is small, then the flame is extinguished
--- the solution {\it quenches}. If, on the other hand, the
power of decay is smaller than three or initial datum is large,
then quenching does not happen, and the burning region spreads
linearly in time. This extends results of Aronson-Weinberger for
the no-flow case. We also consider shear flows with large
amplitude and show that if the reaction power-law decay is larger
than three and the flow has only small plateaux (connected domains
where it is constant), then any compactly supported initial datum
is quenched when the flow amplitude is large enough (which is not
true if the power is smaller than three or in the presence of a
large plateau). This extends results of Constantin-Kiselev-Ryzhik
for combustion with ignition temperature cutoff. Our work carries
over to the case $\Omega = \mathbb{R}^n\times\mathbb{T}^m$, when
the critical power is $1+\tfrac 2n$, as well as to certain
non-periodic flows.
\end{abstract}

\section{Introduction} \lb{S1}

We study the reaction-diffusion-advection equation
\begin{align}
T_t + u\cdot\nabla T &= \Delta T + M f(T) \lb{1.1}
\\ T(0,x) & = T_0(x)\ge 0 \notag
\end{align}
in $\Omega\subseteq\bbR^n$, which models flame propagation in a
fluid \cite{Ber-Lar-Lions} advected by a prescribed flow $u(x)$.
Here $T$ is the normalized temperature that takes values in
$[0,1]$ and $f:[0,1]\to\bbR^+_0$ with $f(0)=f(1)=0$ is the
non-linear reaction term, with coupling $M>0$. There is a vast
mathematical and physical literature on the subject and we refer
to recent reviews \cite{Berrev,Xin2} for an extensive
bibliography. In the present paper we will mainly focus on the
question of {\it quenching} (extinction) of the flame
\begin{equation} \lb{1.1b}
\lim_{t\to\infty} \|T(t,\cdot)\|_\infty =0
\end{equation}
(in which case we say that $T$ {\it quenches}), or its absence.
This means that we will assume the spatial domain to be unbounded
and the initial datum compactly supported (fast enough decay at
infinity would be sufficient). That is, the fluid will be
initially ``hot'' in a finite (but possibly large) central region
and ``cold'' at infinity. Our main interest is in the study of
situations when quenching depends on the size of (the support of)
the initial datum (Theorems~\ref{T.1.3} and \ref{T.3.1},
Corollary~\ref{C.2.3}), or when it results from strong fluid
motion (Theorem~\ref{T.1.2}).

We mainly want to consider {\it combustion-type} reaction terms
with $f'(0)=0$. However, unlike most previous works studying
quenching in reaction-diffusion models, we will not assume $f$ to
have an {\it ignition temperature cutoff}, that is, we will not
require the existence of $\tht_0>0$ such that $f(T)=0$ for
$T\in[0,\tht_0]$. Such an assumption simplifies the proof of
quenching to showing the existence of a time $t_0$ at which $T$ is
below the ignition temperature $\tht_0$, uniformly in space. Then
the maximum principle shows that this will remain the case for all
later times and we are left with a linear equation after $t_0$.
Quenching is now provided by the diffusion term $\Delta$.

Dispensing with this assumption allows us to treat the important
case of {\it Arrhenius reaction term} $f(T)\equiv e^{-c/T}$, and
more generally, our quenching results will hold when $f(T)\le
cT^{p}$ for certain $p>1$. Without the ignition temperature cutoff
the equation will never become linear but can be close to it when
$T$ is small. The idea is that if at low temperatures the reaction
is weak (i.e., if $p$ is large enough), then the decay of
temperature caused by diffusion may still be sufficient to ensure
quenching. Hence we will consider the non-linear equation as a
perturbation of its linear counterpart
\begin{equation} \lb{1.5}
\Phi_t + u\cdot\nabla \Phi = \Delta \Phi.
\end{equation}
Then we will apply a lemma of Meier (Lemma \ref{L.2.1} below) to
show that the solution of the latter can be used to estimate that
of the former. (There is the obvious estimate $T(t,x,y)\le e^{ct}
\Phi(t,x,y)$ with $c\equiv M\|f(T)/T\|_\infty$, following from the
maximum principle, but it is insufficient for our purposes.) To do
all this we will need good estimates on the decay of the solutions
of \eqref{1.5}, which enter into Lemma \ref{L.2.1}.

The first work studying the extinction and propagation of flames
in the case of combustion non-linearity with ignition temperature
cutoff was the paper \cite{Kanel} by Kanel', who considered
equation \eqref{1.1} in one spatial dimension and with no
advection. He showed that if the initial condition is
$T_0(x)\equiv \chi_{[-L,L]}(x)$, then there are two length scales
$L_0$, $L_1$ such that flame extinction/propagation happens when
$L$ is smaller than $L_0$/larger than $L_1$. That is,
\begin{align*}
& T(t,x)\to 0 \text{ as $t\to \infty$ uniformly in $x\in\bbR$ if
$L<L_0$,}
\\ & T(t,x)\to 1 \text{ as $t\to \infty$ for all $x\in\bbR$ if
$L>L_1$.}
\end{align*}
Both length scales are of the order of the {\it laminar front
width} $\ell_c\equiv M^{-1/2}$. But quenching often operates on
larger scales, especially in the presence of strong fluid motion
(see Theorem \ref{T.1.2}).

Kanel's result was generalized by Roquejoffre \cite{Roq} to the
case of {\it shear flows} $u(x,y)=(u(y),0)$ in a cylindrical
domain $\bbR\times D$ with $D\subset\bbR^{m-1}$ bounded and
Neumann boundary conditions at $\partial D$. The length scales
$L_0$, $L_1$ then also depend on $u$. Xin \cite{Xin} extended the
propagation part of Kanel's result to smooth periodic flows on
$\bbR\times[0,h]^m$ with periodic boundary conditions. The
following theorem is an extension of these results to the case of
combustion without ignition temperature cutoff, when $u$ is a
periodic flow on $\bbR\times[0,h]^m$. It identifies the critical
exponent $p^*\equiv 3$ such that the above extinction--propagation
dichotomy picture is valid when $p>p^*$ and $f(T)\le cT^p$ close
to $T=0$, whereas if $p<p^*$ and $f(T)\ge cT^p$ close to $T=0$,
then no non-trivial non-negative solution of \eqref{1.1} satisfies
\eqref{1.1b}.

\begin{theorem} \lb{T.1.3}
Consider \eqref{1.1} in $\Omega\equiv\bbR\times[0,h]^m$ with
periodic boundary conditions. Let $u(x,y)$ be a smooth, periodic,
divergence-free, mean-zero flow on $\Omega$, and let $f$ be
Lipshitz with $f(0)=f(1)=0$ and $f(T)>0$ for $T\in(0,1)$. Let
$c,\eta,\theta>0$ and assume $0\le T_0\le 1$.
\begin{SL}
\item[{\rm{(i)}}] There are $0<\gamma_1<\gamma_2<\infty$,
independent of $\eta$, and $L_1(\eta)<\infty$ such that if
$T_0(x,y)\ge \eta\chi_{[-L_1,L_1]}(x)$ is compactly supported, the
solution of \eqref{1.1} satisfies
\begin{align}
\lim_{t\to\infty} \inf_{|x|\le\gamma_1 t} T(t,x,y) & = 1,
\lb{1.5a}
\\ \lim_{t\to\infty} \sup_{|x|\ge\gamma_2 t} T(t,x,y) & = 0.
\lb{1.5b}
\end{align}
\item[{\rm{(ii)}}]  If $p>3$ and $f(T)\le cT^{p}$ for $0\le
T\le\theta$, then there is $\eps>0$ such that if
$\|T_0\|_1\le\eps$, then the solution of \eqref{1.1} quenches.
\item[{\rm{(iii)}}] If $p<3$ and $f(T)\ge cT^{p}$ for $0\le
T\le\theta$, then the solution of \eqref{1.1} quenches only if
$T_0\equiv 0$. Moreover, there are $0<\gamma_1<\gamma_2<\infty$
such that any solution $T$ with compactly supported $T_0\not\equiv
0$ satisfies \eqref{1.5a}, \eqref{1.5b}.
\end{SL}
\end{theorem}

{\it Remarks.} 1. Part (i) is essentially a result of Xin
\cite{Xin} and we only include it for the sake of completeness.
Part (iii) for $p=1$ and $u$ a shear flow was proved by
Roquejoffre \cite{Roq}.
\smallskip

2. We only need to assume $u$ smooth and divergence-free in part
(i) (we assume $u$ to be smooth when $\Omega$ is viewed as
$\bbR\times(h\bbT)^m$; similarly in Theorem \ref{T.1.2} where $u$
is $C^1$). The mean-zero assumption is not essential, as any
periodic flow is mean-zero in a suitable moving frame.
\smallskip

3. Parts (ii) and (iii) extend to the case
$\Omega\equiv\bbR^n\times [0,h]^m$, with the critical exponent
being $p^*\equiv 1+\tfrac 2n$, as follows from Theorem
\ref{T.3.1}.
\smallskip

In Theorem \ref{T.1.3} quenching results from smallness of the
initial datum, thanks to which $T$ quickly becomes small enough so
that the effects of reaction are weak. On the other hand, large
initial flames can be extinguished by a strong wind.
Constantin-Kiselev-Ryzhik \cite{CKR} studied quenching by large
amplitude shear flows, and considered the problem
\begin{align}
T_t + Au(y) T_x & = \Delta T + M f(T) \lb{1.2}
\\ T(0,x,y) & = T_0(x,y)\ge 0 \notag
\end{align}
on the strip $\bbR\times[0,h]$ with periodic boundary conditions
and flow amplitude $A$. Their interest was in identifying flow
profiles $u$ such that quenching happens for any compactly
supported $T_0$ when $A$ is large enough. They made the following
definition.

\begin{definition} \lb{D.1.1} We say that the profile $u$ is
{\it quenching} if for any compactly supported $T_0(x,y)$, there
exists $A_0$ such that for all $|A|\ge A_0$ the solution of
\eqref{1.2} quenches.
\end{definition}

Of course, whether $u$ is quenching depends on $f$ and $M$. Under
the ignition temperature cutoff assumption on $f$ it is proved in
\cite{CKR} that if a $C^\infty$ profile $u$ has no {\it plateaux}
(intervals on which $u$ is constant) or has only one small
plateau, then it is quenching. On the other hand if $u$ has a
large enough plateau, then it is not quenching. Both these plateau
sizes depend on $f$ and $M$.

Kiselev-Zlato\v s \cite{KZ} later obtained a sharp result in this
direction by showing that there is a critical length $\ell_0(f,M)$
such that $u\in C^1(h\bbT)$ is quenching when all its plateaux are
shorter than $\ell_0$ and it is not quenching when at least one
plateau is longer than $\ell_0$. They also provided estimates on
the minimal {\it quenching amplitude} $A_0$ as a function of the
size of the support of $T_0$ and studied the dependence of this
relation on the (large and small period) scaling of the flow
profile in $y$. All their results agree with previously obtained
numerical experiments (see, e.g., \cite{VCKRR}). Finally,
quenching by large amplitude cellular flows was recently studied
by Fannjiang-Kiselev-Ryzhik \cite{FKR}.

The following theorem is an extension of the results in
\cite{CKR,KZ} to the case of combustion without ignition
temperature cutoff, when $u$ is a shear flow on
$\bbR\times[0,h]^m$ (in which case plateaux of $u$ are connected
sets in $(h\bbT)^m$ on which $u$ is constant). It again identifies
the critical exponent $p^*\equiv 3$ such that the above
quenching--non-quenching dichotomy picture is valid when $p>p^*$
and $f(T)\le cT^p$ close to $T=0$, whereas if $p<p^*$ and $f(T)\ge
cT^p$ close to $T=0$, then quenching never happens.

\begin{theorem} \lb{T.1.2}
Consider \eqref{1.2} on $\Omega\equiv \bbR\times[0,h]^m$ with
periodic boundary conditions. Let $u(x,y)=(u(y),0)$ be a $C^1$
shear flow profile on $\Omega$ and let $0\le f\not\equiv 0$ be
Lipshitz with $f(0)=f(1)=0$. Let $c,\theta>0$.
\begin{SL}
\item[{\rm{(i)}}] If $u$ has at least one large enough (depending
on $f,M$) plateau, then $u$ is not quenching.
\item[{\rm{(ii)}}]
If $p>3$ and $f(T)\le cT^{p}$ for $0\le T\le\theta$, and if $u$
has none or only small enough (depending on $f,M$) plateaux, then
$u$ is quenching.
\item[{\rm{(iii)}}] If $p<3$ and $f(T)\ge
cT^{p}$ for $0\le T\le\theta$, then $u$ is not quenching.
\end{SL}
\end{theorem}

{\it Remarks.} 1. Part (i) for $m=1$ is a result of
Constantin-Kiselev-Ryzhik \cite{CKR}.
\smallskip

2. Large/small enough plateau in (i)/(ii) means one
containing/contained in a large/small enough ball in $(h\bbT)^m$.
The change of variables $\til T(t,x,y)\equiv
T(M^{-1}t,M^{-1/2}x,M^{-1/2}y)$ shows that bounds on the sizes of
both balls (upper on the large one and lower on the small one) are
of the order of the laminar front width $\ell_c\equiv M^{-1/2}$
for any fixed $f$.
\smallskip

3. This result holds with Neumann boundary conditions as well. It
also generalizes to shear flows $u(x,y)=(u(y),0)$ on
$\bbR^n\times[0,h]^m$. The critical exponent is then $p^*\equiv
1+\tfrac 2n$.
\smallskip

4. If $u$ is mean-zero, then \eqref{1.5a},\eqref{1.5b} hold in
(i) and (iii) (in (i) by extension of an argument from \cite{CKR},
in (iii) by Theorem \ref{T.3.1}).
\smallskip

The second group of papers addressing problems related to ours
study the semi-linear heat equation
\begin{equation} \lb{1.4}
T_t + u\cdot \nabla T = \Delta T + T^p
\end{equation}
with $p>1$ on $\bbR^n$, and the first of them was the work of
Fujita \cite{Fuj}. In the case $u\equiv 0$ he showed that if $p>
p^*\equiv 1+\tfrac 2n$, then there are global positive solutions
to \eqref{1.4}, whereas if $1<p<p^*$, then all non-trivial
non-negative solutions blow up in finite time. The critical case
$p=p^*$ was shown to belong to the blowup regime by Hayakawa
\cite{Hay}.

Bandle-Levine \cite{BL} extended Fujita's result to divergence
free flows with $x^{-1}$ decay at infinity, and the existence of a
critical exponent $p^*$ for any flow was proved by Meier
\cite{Me}. In both of these works the Hayakawa case $p=p^*$ is
left open. Several authors have studied the problem on conical or
general sectorial domains, or with additional potential or
non-linear terms in \eqref{1.4}. We refer to the reviews by Levine
\cite{Lev} and Deng-Levine \cite{DL} for more details and
bibliography. In this direction we prove Corollary \ref{C.2.3}
which extends Fujita's theorem to more general classes of flows,
periodic in particular, and is a direct application of lemmas by
Meier \cite{Me} and Norris \cite{No}. It shows that in $\bbR^n$,
the critical exponent for these flows is again $p^*\equiv 1+\tfrac
2n$.

The rest of the paper is organized as follows. In Section \ref{S2}
we state the abovementioned lemmas of Meier and Norris, and their
consequence, Corollary \ref{C.2.3}. In Section \ref{S3} we prove a
general extinction--propagation result (Theorem \ref{T.3.1}), as
well as Theorems \ref{T.1.3} and \ref{T.1.2}.

For the sake of simplicity of notation, in what follows we will be
studying the equation
\[
T_t  = \Delta T + u\cdot\nabla T + f(T)
\]
in $\bbR^n\times\bbT^m$ instead of \eqref{1.1} in
$\bbR^n\times[0,h]^m$ with periodic boundary conditions. This is
no loss as one can be obtained from the other by a change of
variables. Indeed --- if $T$ satisfies \eqref{1.1} in
$\bbR^n\times[0,h]^m$, then $\til T(t,x,y)\equiv T(h^2t,hx,hy)$
satisfies
\[
\til T_t = \Delta\til T + v\cdot\nabla T + g(T)
\]
in $\bbR^n\times\bbT^m$, with $v(x,y)\equiv -h u(hx,hy)$ and
$g(T)\equiv h^2 M f(T)$.

The author would like to thank Alexander Kiselev, James Norris,
Yehuda Pinchover, and Peter Pol\' a\v cik for valuable
communications.

\section{Lemmas of Meier and Norris} \lb{S2}

We now state a lemma of Meier \cite{Me} which enables one to treat
certain reaction-diffusion non-linear PDE's as perturbations of
associated linear equations when one is interested in qualitative
phenomena like extinction and blowup. We state it in the form we
will need here and provide the proof for later reference.

We let $\Omega\subseteq \bbR^n$ be a domain with a piecewise
smooth (possibly empty) boundary $\partial\Omega$. We assume that
$u:\Omega\to\bbR^n$ and $f:\bbR^+_0\to\bbR^+_0$ with $f(0)=0$ are
bounded and $f$ is Lipshitz. We let $T(t,x)$, $\Phi(t,x)$ be the
solutions of
\begin{align}
T_t & = \Delta T + u\cdot\nabla T + f(T) \label{2.1}
\\ \Phi_t & = \Delta \Phi + u\cdot\nabla\Phi     \label{2.2}
\end{align}
on $\Omega$ with Dirichlet, Neumann, or periodic boundary
conditions at $\partial\Omega$, and initial conditions
$T_0(x),\Phi_0(x)\ge 0$ (hence, by the maximum principle,
$T,\Phi\ge 0$). In this section $\|\cdot\|$ stands for
$\|\cdot\|_\infty$.

\begin{lemma}[Meier] \lb{L.2.1} Consider $T$ and $\Phi$ as above and let $c,\al>0$.
\begin{SL}
\item[{\rm{(i)}}] If $f(T)\le cT^{1+\alpha}$ and $I\equiv \int_0^\infty \|\Phi(t,\cdot)\|^\alpha
\,dt$,
then for $0\le\del_0<(c\alpha I)^{-1/\alpha}$ and $T_0(x)\equiv
\del_0\Phi_0(x)$ the solution $T$ quenches.
\item[{\rm{(ii)}}] If $f(T)\ge cT^{1+\alpha}$ and $J\equiv
\sup_{t,x} t \Phi(t,x)^\alpha$,
then for $\del_0>(c\alpha J)^{-1/\alpha}$ and $T_0(x)\equiv
\del_0\Phi_0(x)$ the solution $T$ blows up in finite time.
\end{SL}
\end{lemma}


{\it Remarks.} 1. A more general form is valid with $f(T)$
replaced by $h(t)f(T)$ where $h$ is non-negative and continuous
(see \cite{Me}). In this case $I\equiv \int_0^\infty
h(t)\|\Phi(t,\cdot)\|^\alpha dt$ and $J\equiv \sup_{t,x}
\Phi(t,x)^\alpha \int_0^t h(s)ds$. If $h\in L^1(\bbR^+)$, we also
need $\|\Phi(t,\cdot)\|\to 0$ in (i). Meier only considers the
non-linear term $h(t) T^{1+\alpha}$ but the general case is
identical.
\smallskip

2. In our applications $\Omega$ is unbounded and decay of $\Phi$
will be provided by the diffusion term in \eqref{2.1}.
\smallskip

3. We note that one can replace $u\cdot\nabla T$ by a $C^1$
function $g(t,x,\nabla T)$ as long as $g(t,s,0)=0$ and
$g(t,x,sv)\ge sg(t,x,v)$ for any $v\in\bbR^n$ and $s\ge s_0$, in
which case we also need $\del_0\le (s_0^\alpha + c\alpha
I)^{-1/\alpha}$ in (i) and $\del_0\ge s_0$ in (ii). Interestingly
enough, if instead $g(t,s,0)=0$ and $g(t,x,sv)\le sg(t,x,v)$ for
any $v\in\bbR^n$ and $s\ge s_0$, and $(c\alpha
I)^{-1/\alpha}>s_0$, then the conclusion of (i) is still valid
--- by first obtaining it as below for $s_0\le\del_0<(c\alpha
I)^{-1/\alpha}$ and then for all smaller $\del_0$ by comparison
theorems (see, e.g., \cite[Chapter 10]{Sm}).

\begin{proof}
(i) We can assume $I<\infty$, otherwise there is nothing to prove.
Let $\del(t)$ with $\del(0)\equiv \del_0$ solve
\[
\del'(t) = c\|\Phi(t,\cdot)\|^\alpha \del(t)^{1+\alpha}
\]
so that
\[
\del(t)=\bigg( \del_0^{-\alpha} - c\alpha\int_0^t
\|\Phi(s,\cdot)\|^\alpha\,ds \bigg)^{-1/\alpha}.
\]
If $\del_0^{-\alpha}>c\alpha I$, then $\del(t)$ exists and is
bounded for all $t\in\bbR^+_0$. Now define $\til T(t,x)\equiv
\del(t)\Phi(t,x)$. Then
\[
\til T_t = \Delta \til T + u\cdot\nabla \til T +
c\del^{1+\alpha}\Phi\|\Phi\|^\alpha,
\]
so $\til T$ is a supersolution of \eqref{2.1} with $\til T_0=T_0$,
and we have $\til T\ge T$. Since $\|\Phi(t,\cdot)\|$ is
non-increasing by the maximum principle, $I<\infty$ gives
$\|\Phi(t,\cdot)\|\to 0$. Hence $\|\til T(t,\cdot)\|\to 0$ and the
same is true for $T$.

(ii) Let $w(t,\phi)$ solve
\[
\frac{\partial w}{\partial t}=cw^{1+\alpha}
\]
with $w(0,\phi)\equiv \phi\ge 0$ and define $\til T(t,x) \equiv
w(t,\del_0\Phi(t,x))$ so that
\[
\til T_t = \Delta \til T - \frac{\partial^2 w}{\partial\phi^2}
\del_0^2 |\nabla\Phi|^2 + u\cdot\nabla \til T + c\til
T^{1+\alpha}.
\]
Now
\[
w(t,\phi) = \big( \phi^{-\alpha} - c\alpha t \big)^{-1/\alpha},
\]
so $\tfrac{\partial^2 w}{\partial\phi^2}\ge 0$ and hence $\til T$
is a subsolution of \eqref{2.1}. Since $\til T_0=T_0$, we have
$\til T\le T$. Finally, blow-up of $\til T$ (and of $T$) is
guaranteed by the existence of $t,x$ such that
$(\del_0\Phi(t,x))^{-\alpha}\le c\alpha t$, which follows from
$\del_0>(c\alpha J)^{-1/\alpha}$.
\end{proof}

To apply Lemma \ref{L.2.1} we need to obtain good large-time
asymptotic estimates of heat kernels corresponding to certain
linear equations. One such result is the following lemma of Norris
\cite{No} (for a proof see Theorem~1.1 in \cite{No}). We start
with

\begin{definition} \lb{D.2.1a}
A function $u:\bbR^n\to\bbR^n$ is of {\it type (N)} if $u(x)\equiv
\divg\beta(x) + ({\rm Id}-\beta(x))\nabla\log\mu(x) + \bar
b/\mu(x)$ with $\beta$ a bounded, differentiable, and
antisymmetric $n\times n$ matrix, $\mu$ positive, differentiable,
and bounded away from $0$ and $\infty$, and  $\bar b\in\bbR^n$ a
constant vector. If $\bar b\neq 0$, we also require the existence
of a bounded, differentiable vector field $\xi$ such that
$\divg(\mu\xi)+\mu\equiv 1$. By the discussion on p.~168 of
\cite{No}, this includes all $u$ periodic with bounded $\divg u$.
\end{definition}

{\it Remark.} Theorems \ref{T.1.3} and \ref{T.1.2} involve
periodic divergence-free $u$. Such functions of type (N) can be
written as $u(x)\equiv \divg\beta(x) + \bar b$ (see
\cite[p.~168]{No}), and so the {\it effective drift} $\bar b$ is
just the mean of $u$.

\begin{lemma}[Norris] \lb{L.2.2}
If $u:\bbR^n\to\bbR^n$ is of type (N), then there is $C<\infty$
such that for any $x,y\in\bbR^n$ and $t>0$, the heat kernel
$k(t,x,y)$ of \eqref{2.2} in $\bbR^n$ satisfies
\begin{equation} \label{2.4}
C^{-1}t^{-n/2} e^{-C|x-y|^2/t} \le k(t,x,y+\bar bt) \le Ct^{-n/2}
e^{-|x-y|^2/Ct}.
\end{equation}
\end{lemma}

{\it Remark.} Of course, $k$ is such that
\[
\Phi(t,x) = \int_{\bbR^n} k(t,x,y) \Phi_0(y) dy.
\]
\smallskip

As an immediate application of Lemmas \ref{L.2.1} and \ref{L.2.2}
we obtain a generalization of a result of Fujita \cite{Fuj}.

\begin{corollary} \lb{C.2.3}
Let $u:\bbR^n\to\bbR^n$ be $C^1$ and of type (N), and consider
\begin{equation} \label{2.5}
T_t = \Delta T + u\cdot\nabla T + T^{1+\alpha}
\end{equation}
in $\bbR^n$.
\begin{SL}
\item[{\rm{(i)}}] If $\alpha>\tfrac 2n$, then there are global
positive solutions of \eqref{2.5} that quench. \item[{\rm{(ii)}}]
If $0<\alpha<\tfrac 2n$, then all non-trivial non-negative
solutions of \eqref{2.5} blow up in finite time.
\end{SL}
\end{corollary}

{\it Remarks.} 1. Fujita proved this for $u\equiv 0$. In that case
the conclusion of (ii) also holds when $\alpha = \tfrac 2n$
\cite{Hay}. Our result (i) is slightly stronger in that quenching
is provided by small enough $L^1$ and $L^\infty$ norms of $T_0$,
with no additional conditions on its decay.
\smallskip

2. Since both Lemmas \ref{L.2.1} and \ref{L.2.2} hold when
$\Delta$ is replaced by a uniformly elliptic operator
$\sum_{i,j=1}^n a_{i,j}(x) \tfrac{\partial^2}{\partial x_i\partial
x_{j}}$ with bounded differentiable $a_{i,j}(x)$, so does this
corollary. The same is true for Theorem \ref{T.3.1}.

\begin{proof}
(i) Let $\Phi$ be the solution of
\begin{equation} \label{2.6}
\Phi_t = \Delta \Phi + u\cdot\nabla \Phi
\end{equation}
in $\bbR^n$ with initial condition $0<\Phi_0\in L^1(\bbR^n)\cap
L^\infty(\bbR^n)$. Then by \eqref{2.4} and the maximum principle,
\[
\|\Phi(t,\cdot)\| \le \min\{ \|\Phi_0\|,
\|\Phi_0\|_1\|k(t,\cdot,\cdot)\| \} \le \min\{ \|\Phi_0\|,
Ct^{-n/2}\|\Phi_0\|_1  \}.
\]
Hence $I\equiv \int_0^\infty \|\Phi(t,\cdot)\|^\alpha\, dt<\infty$
and Lemma \ref{L.2.1}(i) gives the result.

(ii) Assume $T_0(x_0)>0$ and $\Phi_0\equiv T_0\ge 0$. Then by
\eqref{2.4},
\[
\Phi(t,x_0-\bar bt)= \int_{\bbR^n} k(t,x_0-\bar bt,y) \Phi_0(y)\,
dy \ge C^{-1}t^{-n/2}e^{-C}D
\]
for $t\ge 1$ and $D\equiv \int_{B(x_0,1)} \Phi_0(y)\, dy>0$ (here
$B(x_0,1)$ is the ball in $\bbR^n$ with center $x_0$ and radius
1). But then $J\equiv\infty$ in Lemma~\ref{L.2.1}(ii), so $T$
blows up in finite time.
\end{proof}

\section{Proofs of the main results} \lb{S3}

We now proceed to prove Theorems \ref{T.1.3} and \ref{T.1.2}. We
will start with a general result in the domain
$\bbR^n\times\bbT^m$, which is related to Corollary \ref{C.2.3}.
We will assume $u$ to be $C^1$ and $f:[0,1]\to\bbR^+_0$ to be
Lipshitz with $f(0)=f(1)=0$.

\begin{theorem} \lb{T.3.1}
Consider \eqref{2.1} on $\bbR^n\times\bbT^m$ with $n\ge 1$ and
$m\ge 0$ and let $u:\bbR^n\times\bbT^m\to\bbR^{n+m}$ be of type
(N). Let $c,\theta>0$ and $0\le T_0\le 1$.
\begin{SL}
\item[{\rm{(i)}}] If $\alpha>\tfrac 2n$ and $f(T)\le
cT^{1+\alpha}$ for $0\le T\le\theta$, then there is $\eps>0$ such
that if $\|T_0\|_1\le\eps$, then the solution of \eqref{2.1}
quenches. \item[{\rm{(ii)}}] If $\alpha<\tfrac 2n$ and $f(T)\ge
cT^{1+\alpha}$ for $0\le T\le\theta$, then the solution of
\eqref{2.1} quenches only if $T_0\equiv 0$. If also $f(T)>0$ for
$T\in (0,1)$, then there are $0<\gamma_1<\gamma_2<\infty$ such
that any solution $T$ with compactly supported $T_0\not\equiv 0$
satisfies
\begin{align}
\lim_{t\to\infty} \inf_{|x|\le\gamma_1 t} T(t,x-\bar bt) & = 1,
\lb{3.0a}
\\ \lim_{t\to\infty} \sup_{|x|\ge\gamma_2 t} T(t,x-\bar bt) & = 0,  \lb{3.0b}
\end{align}
with $\bar b$ from Definition \ref{D.2.1a}.
\end{SL}
\end{theorem}

{\it Remarks.} 1. If $f$ is as in (ii), the theorem says that no
flame can be extinguished, even in the presence of strong
advection of type (N).
\smallskip

2. For $u\equiv 0$ in $\bbR^n$ this was proved by
Aronson-Weinberger \cite{AW}, using results from \cite{Fuj}.

\begin{proof}
Let $q(t,x,y)$ be the heat kernel for $\Delta+u\cdot\nabla$ in
$\bbR^n\times\bbT^m$ and $k(t,x,y)$ the one in $\bbR^{n+m}$ (with
$u$ periodically continued in the last $m$ coordinates). Then
\[
q(t,x,y)=\sum_{j\in\bbZ^m} k(t,x,y+(0,j))
\]
where $(0,j)\in \bbR^{n+m}$. Hence from \eqref{2.4} with $n+m$ in
place of $n$ we obtain for all $x,y\in \bbR^n\times\bbT^m$ and
$t\ge t_0>0$
\begin{equation} \lb{3.1}
C^{-1}t^{-n/2} e^{-C|x-y|_*^2/t} \le q(t,x,y+\bar bt) \le
Ct^{-n/2} e^{-|x-y|_*^2/Ct}
\end{equation}
with some new $C=C(u,t_0)<\infty$. Here $\bar b\in\bbR^{n+m}$ and
 $|x |_*$ denotes the norm of the $\bbR^n$ component of $x$.
In particular,
\[
C^{-1}t^{-n/2} \le \|q(t,\cdot,\cdot)\|_\infty \le Ct^{-n/2}
\]
for $t\ge t_0$.

(i) By changing $c$ we can assume $f(T)\le cT^{1+\alpha}$ for all
$T$. Let $\Phi$ satisfy \eqref{2.2} in $\bbR^n\times\bbT^m$ with
$\Phi_0\equiv T_0$. Then we have $\|\Phi(t,\cdot)\|_\infty \le 1$,
and for $t\ge t_0$
\[
\|\Phi(t,\cdot)\|_\infty \le
\|\Phi_0\|_1\|q(t,\cdot,\cdot)\|_\infty \le C\eps t^{-n/2}
\]
Hence we get
\[
I \equiv \int_0^\infty \|\Phi(t,\cdot)\|_\infty^\alpha \,dt \le
t_0 + \int_{t_0}^\infty (C\eps t^{-n/2})^\alpha\, dt <\frac
1{c\alpha}
\]
if $t_0$ and then $\eps$ are chosen small enough. Lemma
\ref{L.2.1}(i) with $\del_0=1$ then gives the result.

(ii) It is obviously sufficient to consider $\alpha>0$. Assume
$\Phi_0\equiv T_0\not\equiv 0$ so that $0<\Phi(t,x)\le T(t,x)$ for
all $x$ and $t>0$. If \eqref{1.1b} were true, we would have
$T(t,x)\le\theta$ and thus $f(T(t,x))\ge cT(t,x)^{1+\alpha}$ for
all $x$ and $t\ge t_1>0$. Let $\til\Phi(t,x)\equiv\Phi(t+t_1,x)$
and let $\til T$ be the solution of \eqref{2.1} with $\til
T(0,x)\equiv\til\Phi(0,x)$, so that $\til T(t-t_1,x)\le
T(t,x)\le\theta$ for $t\ge t_1$. Then by \eqref{3.1},
\[
\til\Phi(t,-\bar bt)= \int_{\bbR^n\times\bbT^m} q(t,-\bar bt,y)
\til\Phi_0(y)\, dy \ge C^{-1}t^{-n/2}e^{-C}D
\]
for $t\ge 1$ and $D\equiv \int_{B(0,1)\times\bbT^m}
\til\Phi_0(y)\, dy>0$. But then $J\equiv\infty$ in Lemma
\ref{L.2.1}(ii), so $\til T$ blows up in finite time, a
contradiction.

Let us now prove \eqref{3.0a} For the sake of transparency we will
assume $m=0$. The proof in the general case is identical, with all
domains $D$ replaced by $D\times\bbT^m$. Let us also assume $\bar
b=0$. Otherwise one can make the change of variables $\til
T(t,x)=T(t,x-\bar bt)$ which gives $\til u(x)\equiv u(x)+\bar b$
in \eqref{2.1} for $\til T$ and hence $\til k(t,x,y)\equiv
k(t,x-\bar bt,y)=k(t,x-\bar bt,y-\bar bt+\bar bt)$ in \eqref{2.4}
for $\til u$. Therefore one has $\til b=0$ in this problem and
proving the claim for $\til T$ immediately gives it for $T$ as
well.

Since $f(T)>0$ for $T\in (0,1)$, we can change $c>0$ so that we
can take $\theta\equiv \tfrac 12$. If $0\le T\not\equiv 0$, we
have $\|T_0\chi_{B(x_0,1)}\|_1=Ce^{9C}\eps$ for some $x_0$ and
$\eps>0$. Assume for simplicity that $x_0=0$ (all the following
estimates will be uniform in $x_0$). If we let $\Phi$ satisfy
\eqref{2.2} with $\Phi_0\equiv T_0$, we have for $\tau^2\ge 1$
\begin{equation} \lb{3.1a}
\Phi(\tau^2,x)\ge \bar\Phi_0(x)\equiv \eps \tau^{-n}
\chi_{B(0,2\tau)}(x)
\end{equation}
by \eqref{3.1} for $t_0\equiv 1$. Obviously if we let $\bar
T,\bar\Phi$ satisfy \eqref{2.1},\eqref{2.2} with initial data
$\bar T_0\equiv\bar \Phi_0$, we have $T(\tau^2+t,x)\ge \bar
T(t,x)$ by comparison theorems, so it is sufficient to prove the
claim for $\bar T$.

Next let
\[
\til T(t,x)\equiv w(t,\bar\Phi(t,x))=
\big(\bar\Phi(t,x)^{-\alpha}-c\alpha t \big)^{-1/\alpha}.
\]
We obviously have $\bar\Phi \le \eps\tau^{-n}$ and so $\til
T(t,x)\le (\eps^{-\alpha}\tau^{\alpha n}-c\alpha t )^{-1/\alpha}$.
Hence up to time $t^2\equiv (\eps^{-\alpha}\tau^{\alpha
n}-2^\alpha)/c\alpha$ we have $\til T\le \tfrac 12$ (if
$\tau>(2\eps)^{1/n}$). So by the argument in Lemma \ref{L.2.1},
$\til T(t^2,x)\le\bar T(t^2,x)$. Moreover, if $\tau$ is large, we
have $t\le \tau^\omega$ for some $\omega< 1$. Thus we can take
$\tau$ large enough so that
\[
\int_{\bbR^n\smallsetminus B(0,\tau)} C t^{-n} e^{-|y|^2/Ct^2} dy
< \frac{(2\eps)^\alpha}{2\alpha} \tau^{-\alpha n}
\]
(notice that $t^{-n}\sim \tau^{-\alpha n}$). Then we have for any
$x$
\[
\int_{B(x,\tau)} q(t^2,x,y) dy > 1-\frac{(2\eps)^\alpha}{2\alpha}
\tau^{-\alpha n}
\]
by \eqref{3.1} with $\bar b=0$ because the integral over $\bbR^n$
is 1.

Now for $x\in B(0,\tau)$
\[
\bar\Phi(t^2,x)\ge\int_{B(x,\tau)} q(t^2,x,y) \eps \tau^{-n} dy\ge
\bigg( 1-\frac{(2\eps)^\alpha}{2\alpha} \tau^{-\alpha n} \bigg)
\eps \tau^{-n}
\]
and so
\[
\til T(t^2,x)= \big(\bar\Phi(t^2,x)^{-\alpha}-c\alpha t^2
\big)^{-1/\alpha}\ge \frac 14,
\]
using
\[
\bigg( 1-\frac{(2\eps)^\alpha}{2\alpha} \tau^{-\alpha n}
\bigg)^{-\alpha} \le 1+(2\eps)^\alpha \tau^{-\alpha n}
\]
for large $\tau$, and the definition of $t$. Hence for large
enough $\tau$ and $t(\tau)$ as above we have
\[
T(\tau^2+t^2,x)\ge \bar T(t^2,x)\ge \til T(t^2,x) \ge \frac 14
\]
for $x\in B(0,\tau)$. And since $\tau^2+t(\tau)^2$ is continuous
in $\tau$ and belongs to $(\tau^2,2\tau^2)$ when $\tau$ is large,
this implies
\begin{equation} \lb{3.1c}
T(2\tau^2,x)\ge \frac 14 \chi_{B(x_0,\tau)} (x)
\end{equation}
for $x_0=0$ and all $\tau\ge\tau_0$. All the above estimates are
uniform in $x_0$, and so \eqref{3.1c} holds for all $x_0$ and
$\tau\ge\tau_0$, with $\tau_0=\tau_0(\delta)$ depending only on
$\delta\equiv \|T_0\chi_{B(x_0,1)}\|_1$ (and $u,f$, of course).

Let $T_0$ be given and assume $T_0(x)\ge \tfrac 14
\chi_{B(x_0,1)}(x)$ for some $x_0$ (otherwise first pick $t$ so
that $T(t,x)\ge\tfrac 14 \chi_{B(x_0,1)}(x)$ and then reset
$T_0(x)$ to $T(t,x)$).
If $\tau\ge \tau_1\equiv \tau_0(|B(0,1)|/4)$, then by the above,
\[
T(2\tau^2,x)\ge \frac 14 \chi_{B(x_0,2)} (x).
\]
Applying this argument again, with initial datum $T(2\tau^2,x)$
instead of $T_0(x)$, we obtain
\[
T(4\tau^2,x)\ge \frac 14 \chi_{B(x_0,3)} (x)
\]
because $T(2\tau^2,x)\ge\tfrac 14 \chi_{B(x_1,1)}(x)$ for any
$x_1\in B(x_0,1)$. Iteration of this gives us
\[
T(2\tau^2j,x)\ge \frac 14 \chi_{B(x_0,j)} (x).
\]
This holds for any $\tau\in[\tau_1,2\tau_1]$ and it follows that
\begin{equation} \lb{3.1d}
T(\gamma t,x)\ge \frac 14 \chi_{B(0,t)} (x).
\end{equation}
for $\gamma\equiv 16\tau_1^2$ and $t\ge |x_0|$.

The proof of \eqref{3.0a} will be finished by yet another
application of the above argument. Let $\eps>0$ be arbitrary and
let $c_\eps>0$ be such that $f(T)\ge c_\eps T^{1+\alpha}$ for
$0\le T\le 1-\eps$. We will show
\begin{equation} \lb{3.1e}
\liminf_{t\to\infty} \inf_{|x|\le t} T(2\gamma t + s_\eps ,x) \ge
1-2\eps
\end{equation}
for some $s_\eps <\infty$, which will imply \eqref{3.0a} with, for
instance, $\gamma_1\equiv (3\gamma)^{-1}$.

Let $\Phi$ solve \eqref{2.1} with initial condition $\Phi_0\equiv
\tfrac 14 \chi_{B(0,2t)}$ for some $t\ge |x_0|$, and let
\begin{equation} \lb{3.1f}
\til T(s,x)\equiv (\Phi(s,x)^{-\alpha}-c_\eps \alpha
s)^{-1/\alpha}.
\end{equation}
Since $\Phi\le \tfrac 14$, up to time $s_\eps\equiv
(4^\alpha-(1-\eps)^{-\alpha})/c_\eps \alpha$ we have $\til T\le
1-\eps$, and so by the proof of lemma \ref{L.2.1}, \eqref{3.1d},
and comparison theorems, $\til T(s_\eps,x)\le T(2\gamma
t+s_\eps,x)$. For all $x\in B(0,t)$ we have by \eqref{3.1},
\[
\Phi(s_\eps,x)\ge \frac 14 \bigg( 1 - \int_{\bbR^n \smallsetminus
B(0,t)} Cs_\eps^{-n/2}e^{-|y|^2/Cs_\eps} dy \bigg) > \big(
4^\alpha + (1-2\eps)^{-\alpha} - (1-\eps)^{-\alpha}
\big)^{-1/\alpha}
\]
if $t$ is large enough. Plugging this into \eqref{3.1f}, we obtain
$\inf_{|x|\le t} \til T(s_\eps,x)\ge 1-2\eps$ (for any large $t$).
This gives \eqref{3.1e}, and \eqref{3.0a} is proved.

We are left with \eqref{3.0b}. Since $f$ is Lipshitz, there is $d$
such that $f(T)\le dt$. Then by the maximum principle,
\[
T(t,x)\le e^{dt}\Phi(t,x),
\]
with $\Phi$ solving \eqref{2.2} and $\Phi_0\equiv T_0$ compactly
supported. By \eqref{3.1},
\[
\Phi(t,x-\bar bt)\le \int_{\bbR^n\times\bbT^m}
Ct^{-n/2}e^{-|x-y|_*^2/Ct}\Phi_0(y) dy,
\]
which is less than $t^{-n/2}e^{-dt}$ whenever $|x|_*\ge
\sqrt{2Cd\,}\,t$ and $t$ is large. The proof is finished.
\end{proof}

\begin{proof}[Proof of Theorem \ref{T.1.3}]
(i) Eq.~\eqref{1.5a} follows from the same result for combustion
with ignition temperature cutoff $\theta_0<\eta$ \cite{Xin} and
comparison theorems. Eq.~\eqref{1.5b} is just \eqref{3.0b} because
by the remark after Definition~\ref{D.2.1a}, $\bar b=0$.

(ii), (iii) Follow directly from Theorem \ref{T.3.1}(i),(ii) with
$\alpha\equiv  p-1$ and $n\equiv 1$ (since $\bar b=0$).
\end{proof}

Now we turn to the proof of Theorem \ref{T.1.2}. Hence $T$ and
$\Phi$ will be the solutions of
\begin{align}
T_t  & = \Delta T + Au(y)T_x + f(T) \lb{3.2}
\\ \Phi_t & =\Delta \Phi + Au(y)\Phi_x \lb{3.3}
\end{align}
in $\bbR\times\bbT^m$. We will consider the initial condition
\begin{equation} \lb{3.4}
T_0(x,y)\equiv \Phi_0(x,y) \equiv \chi_{[-L,L]}(x),
\end{equation}
since by comparison theorems, $u$ is quenching if and only if for
every $L$ the solution $T$ quenches when $|A|$ is large enough. We
will again use Lemma \ref{L.2.1} but to prove part (ii) we need to
estimate the decay of $\Phi$ without the help of Lemma
\ref{L.2.2}, since the constants in it may not be uniform in $A$.

Instead, we express the solution of \eqref{3.3} in terms of the
Brownian motion. Following \cite{KZ} we obtain $\Phi(t,x,y)  =
\bbE ( \Phi(0,X^x_t,Y^y_t))$ where $\bbE$ is the expectation with
respect to the random process $(X^x_t,Y^y_t)$ starting at $(x,y)$
and satisfying
\begin{align*}
dX^x_t & = \sqrt{2} \, dW^x_t + Au(Y^y_t)dt,
\\ dY^y_t & = \sqrt{2} \, dW^y_t.
\end{align*}
Here $(W^x_t,W^y_t)$ is the normalized Brownian motion on
$\bbR\times\bbT^m$ starting at $(x,y)$. Thus, $Y^y_t=y+\sqrt{2}
(W^y_t-y) = W^y_{2 t}$ and
\[
X^x_t = x+\sqrt{2} (W^x_t-x) + \int_0^t Au(Y^y_s)ds = W^x_{2 t} +
\frac{A}{2}\int_0^{2 t} u(W^y_s)ds.
\]
Then we have by \eqref{3.3}, \eqref{3.4}, and Lemma 7.8 in
\cite{Oks},
\begin{equation} \lb{3.5}
\Phi(t,x,y)  =  \bbP \bigg( W_{2 t}^x + \frac{A}{2}\int_0^{2 t}
u(W_s^y)ds \in[-L, L] \bigg).
\end{equation}

To evaluate this probability we employ Lemmas 2.1 and 2.2 in
\cite{KZ}. There they are proved for $m=1$ but the proof in the
general case is identical.

\begin{lemma}[Kiselev-Zlato\v s] \lb{L.3.2}
If $u\in C^1(\bbT^m)$ and $S\subset (0,\infty)$ is compact, then
\begin{equation*} 
\lim_{\eps\to 0} \sup_{(t,a,y)\in S\times\bbR\times\bbT^m} \bbP
\bigg( \int_0^{2t} u(W_s^y)ds \in [a,a+\eps] \smallsetminus \{
2tu(y) \} \bigg) = 0,
\end{equation*}
while $\bbP ( \int_0^{2t} u(W_s^y)ds = 2tu(y) )$ equals the
probability of $\{W^y_s\}_{s\in[0,2t]}$ staying entirely inside a
plateau of $u$, and is zero unless $y$ is in the interior of a
plateau.
\end{lemma}

In other words, if $0<t_0<t_1<\infty$, then by making $A$ large,
$\Phi(t,x,y)$ can be made as small as we want for $t\in[t_0,t_1]$
and $y$ not in a plateau of $u$, since \eqref{3.5} and
independence of $W^x_{2t}$ and $W^y_s$ imply
\begin{equation} \lb{3.5a}
\Phi(t,x,y)  \le \sup_{a\in\bbR} \bbP \bigg( \int_0^{2 t}
u(W_s^y)ds \in \bigg[a, a+\frac{4L}A \bigg] \bigg).
\end{equation}
This is in line with the intuition that, outside of plateaux of
$u$, strong wind quickly extinguishes the flame by stretching it
and exposing it to diffusion \cite{CKR,KZ}. This takes care of
estimating $\Phi(t,x,y)$ within any finite time interval and for
$y$ not in a plateau. When $y$ is inside a plateau, we also need
the following estimate.

\begin{lemma} \lb{L.3.3}
Let $y\in B(0,\eps)\subseteq\bbT^m$ and $C\equiv 4/\pi$. Then
\begin{equation*} 
\bbP \big( W^y_s\in B(0, \eps) \text{ for all $s\in[0,2t]$} \big)
\le Ce^{-\pi^2 t/4\eps^2}.
\end{equation*}
\end{lemma}

\begin{proof}
This probability is obviously largest when $m=1$ and $y=0$. The
Feynman-Kac formula says that this is
\[
\big( e^{-2tH_\eps} \chi \big) ( 0 ),
\]
where $H_\eps\equiv-\tfrac 12\Delta$ on $[-\eps,\eps]$ with
Dirichlet boundary conditions at $\pm\eps$ and $\chi(y)\equiv
\chi_{[-\eps,\eps]}(y)$. Since on $[-\eps,\eps]$
\[
\chi(y)=\sum_{n=0}^\infty \frac {4(-1)^n}{(2n+1)\pi} \cos \bigg(
\frac {(2n+1)\pi}{2\eps} y \bigg)
\]
and
\[
H_\eps \cos \bigg( \frac {(2n+1)\pi}{2\eps} y \bigg) =
\frac{(2n+1)^2\pi^2}{8\eps^2} \cos \bigg( \frac {(2n+1)\pi}{2\eps}
y \bigg),
\]
we get
\begin{equation*}
\big( e^{-2 tH_\eps} \chi \big) (0)  = \sum_{n=0}^\infty \frac
{4(-1)^n}{(2n+1)\pi} e^{- 2t (2n+1)^2\pi^2/8\eps^2}  \le \frac
4{\pi} e^{-t \pi^2/4\eps^2},
\end{equation*}
since the sum is alternating.
\end{proof}

Finally, large $t$ can be handled by the estimate
\begin{equation} \lb{3.7}
\Phi(t,x,y) \le \sup_{a\in\bbR} \bbP \big( W_{2 t}^x\in[a,a+2L]
\big)\le D t^{-1/2}
\end{equation}
for $D\equiv L\pi^{-1/2}$. The first inequality again follows from
\eqref{3.5} and the independence of $W^x_{2t}$ and $W^y_s$, the
second because the density function of the random variable $W_{2
t}^x$ is $\varphi(z)=(4\pi t)^{-1/2}e^{-|z-x|^2/4t}\le (4\pi
t)^{-1/2}$.

\begin{proof}[Proof of Theorem \ref{T.1.2}]
(i) For $m=1$ this is a result from \cite{CKR}, where a radially
symmetric subsolution of \eqref{2.1} supported on $\bbR\times I$
(for some plateau $I$) is constructed using Bessel functions. When
$m\ge 2$, the same construction can be applied, with an extra
technical difficulty. This stems from the fact that the
fundamental solution of $\Delta T= 0$ in $\bbR^n$ is bounded below
when $n\ge 3$. It can be overcome and the result will follow.

(ii) Again we can change $c$ to get $f(T)\le cT^{1+\alpha}$ for
all $T$, with $\alpha\equiv p-1$.  First assume that $u$ has no
plateaux and let $0<t_0<t_1<\infty$ and
\[
s \equiv\sup_{t\in[t_0,t_1]}
\|\Phi(t,\cdot,\cdot)\|_\infty^{\alpha}.
\]
Then by $\|\Phi(t,\cdot,\cdot)\|_\infty\le 1$ and \eqref{3.7} we
have (with $D\equiv L\pi^{-1/2}$)
\[
I\equiv\int_0^\infty \|\Phi(t,\cdot,\cdot)\|_\infty^{\alpha}\, dt
\le t_0 + s(t_1-t_0) +  \int_{t_1}^\infty D^{\alpha}
t^{-\alpha/2}\, dt.
\]
Now \eqref{3.5a} and Lemma \ref{L.3.2} show that by taking $|A|$
large, on can make $t_0$, $t_1^{-1}$, and $s$ small enough so that
$I<(c\alpha)^{-1}$. The result then follows by taking
$\del_0\equiv 1$ in Lemma \ref{L.2.1}(i).

Let us now assume that $u$ has plateaux, each contained in a ball
of radius less than $\eps<\tfrac \pi{2} d^{-1/2}$, where $d$ is
such that $f(T)\le d T$. Such $d$ exists because $f(0)=0$ and $f$
is Lipschitz. Define $\omega\equiv \pi^2/4\eps^2-d>0$ and let
$0<t_0<t_1<\infty$ be such that $\int_{t_0}^\infty (2Ce^{-\omega
t})^{\alpha}dt<(2c\alpha)^{-1}$ (with $C\equiv 4/\pi$) and
$\int_{t_1}^\infty (e^{d t_0}D t^{-1/2})^{\alpha}dt
<(2c\alpha)^{-1}$. Lemmas \ref{L.3.2} and \ref{L.3.3} show that if
$|A|$ is large enough, then
\begin{equation} \lb{3.8}
\|\Phi(t,\cdot,\cdot)\|_\infty\le 2C e^{-\pi^2t/4\eps^2}
\end{equation}
for $t\in[t_0,t_1]$. By the maximum principle, $T(t,x,y)\le e^{d
t}\Phi(t,x,y)$ and so with $\til T(t,x,y)\equiv T(t+t_0,x,y)$ and
$\til\Phi(t,x,y)\equiv e^{d t_0}\Phi(t+t_0,x,y)$ we have $0\le\til
T(0,x,y)\le\til\Phi(0,x,y)$. By \eqref{3.7} and \eqref{3.8}, we
also have
\begin{align*}
\int_0^\infty \|\til\Phi(t,\cdot,\cdot)\|_\infty^{\alpha} dt & =
\int_{t_0}^\infty \Big( e^{d t_0}\|\Phi(t,\cdot,\cdot)\|_\infty
\Big)^{\alpha} dt
\\ & \le \int_{t_0}^{t_1} \Big( 2Ce^{-\omega t} \Big)^{\alpha} dt +
\int_{t_1}^\infty \Big( e^{d t_0}D t^{-1/2} \Big)^{\alpha} dt
\\ & <(c\alpha)^{-1}.
\end{align*}
Now Lemma \ref{L.2.1}(i) with $\del_0\equiv 1$ gives
$\lim_{t\to\infty} \|\til T(t,\cdot,\cdot)\|_\infty = 0$, and so
the same holds for $T$.

(iii) Follows from Theorem \ref{T.3.1}(ii) with $\alpha\equiv p-1$
and $n\equiv 1$.
\end{proof}


\begin{thebibliography}{99}

\bibitem{AW} D.G.~Aronson and H.F.~Weinberger,
\it Multidimensional nonlinear diffusion arising in population
genetics, \rm Adv. in Math. {\bf 30} (1978), 33--76.

\bibitem{BL} C.~Bandle and H.A.~Levine, \it Fujita type phenomena for
reaction-diffusion equations with convection like terms, \rm
Differential Integral Equations {\bf 7} (1994), 1169--1193.

\bibitem{Berrev} H.~Berestycki, \it The influence of advection on the propagation
of fronts in reaction-diffusion equations, \rm Nonlinear PDEs in
Condensed Matter and Reactive Flows, NATO Science Series C, 569,
H. Berestycki and Y. Pomeau eds, Kluwer, Doordrecht, 2003.

\bibitem{Ber-Lar-Lions} H.~Berestycki, B.~Larrouturou and P.-L.~Lions,
\it Multi-dimensional travelling wave solutions of a flame
propagation model, \rm Arch. Rational Mech. Anal. {\bf 111}
(1990), 33--49.

\bibitem{CKR} P.~Constantin, A.~Kiselev, L.~Ryzhik, \it Quenching
of flames by fluid advection, \rm Comm. Pure Appl. Math. {\bf 54}
(2001), 1320--1342.

\bibitem{DL} K.~Deng and H.A.~Levine, \it The role of critical exponents in
blow-up theorems: the sequel, \rm J. Math. Anal. Appl. {\bf 243}
(2000), 85--126.

\bibitem{FKR} A.~Fannjiang, A.~Kiselev, and L.~Ryzhik, \it
Quenching of reaction by cellular flow, \rm preprint.

\bibitem{Fuj} H.~Fujita, \it On the blowing up of solutions of the Cauchy
problem for $u_{t}=\Delta u+u_{1+\alpha }$, \rm J. Fac. Sci. Univ.
Tokyo Sect. I {\bf 13} (1966) 109--124.

\bibitem{Hay} K.~Hayakawa, \it On nonexistence of global solutions of some
semilinear parabolic differential equations, \rm Proc. Japan Acad.
{\bf 49} (1973), 503--505.

\bibitem{Kanel} Ja.I.~Kanel', \it  Stabilization of the solutions of the
equations of combustion theory with finite initial functions, \rm
Mat. Sb. (N.S.){\bf 65 (107)} 1964 suppl., 398--413.

\bibitem{KZ} A.~Kiselev and A.~Zlato\v s, \it
Quenching of combustion by shear flows, \rm preprint.

\bibitem{Lev} H.A.~Levine, \it The role of critical exponents in blowup
theorems, \rm SIAM Rev. {\bf 32} (1990), 262--288.

\bibitem{Me} P.~Meier, \it On the critical exponent for reaction-diffusion
equations, \rm Arch. Rational Mech. Anal. {\bf 109} (1990),
63--71.

\bibitem{No} J.R.~Norris, \it Long-time behaviour of heat flow: global estimates
and exact asymptotics, \rm Arch. Rational Mech. Anal. {\bf 140}
(1997), 161--195.

\bibitem{Oks} B. Oksendal, \it Stochastic Differential Equations,
\rm Springer-Verlag, Berlin, 1995.

\bibitem{Roq} J.-M.~Roquejoffre, \it Eventual monotonicity and convergence to travelling
fronts for the solutions of parabolic equations in cylinders, \rm
Ann. Inst. H. Poincar\'e Anal. Non Lin\'eaire {\bf 14} (1997),
499--552.

\bibitem{Sm} J.~Smoller, \it Shock Waves and Reaction-Diffusion Equations,
\rm Springer-Verlag, New York, 1994.

\bibitem{VCKRR} N.~Vladimirova, P.~Constantin, A.~Kiselev, O.~Ruchayskiy and L.~Ryzhik,
\it Flame enhancement and quenching in fluid flows, \rm Combustion
Theory and Modelling {\bf 7} (2003), 487--508.

\bibitem{Xin} J.~Xin, \it Existence and nonexistence of travelling waves and reaction-diffusion
front propagation in periodic media, \rm J. Stat. Phys. {\bf 73}
(1993), 893--926.

\bibitem{Xin2} J.~Xin, \it Front propagation in heterogeneous media, \rm SIAM Rev. {\bf 42} (2000),
161--230.

\end{thebibliography}
\end{document}